\documentclass[11pt]{article}

\usepackage[T1]{fontenc}
\usepackage[utf8]{inputenc}
\usepackage{lmodern}
\usepackage{amsmath,amssymb,amsthm,mathrsfs}
\usepackage{mathtools}
\usepackage[margin=1in]{geometry}
\usepackage{microtype}
\usepackage{parskip}
\usepackage{needspace}
\usepackage{xcolor}

\usepackage{enumitem}
\usepackage[hidelinks,bookmarksnumbered=true]{hyperref}

\hypersetup{
  pdftitle={Geometric progressions in distance sets in odd dimensions over prime fields},
  pdfauthor={Thang Pham, Dung The Tran, Boqing Xue}
}
\newtheorem{theorem}{Theorem}[section]
\newtheorem{proposition}[theorem]{Proposition}
\newtheorem{corollary}[theorem]{Corollary}
\newtheorem{lemma}[theorem]{Lemma}
\theoremstyle{remark}

\theoremstyle{definition}
\newtheorem{definition}[theorem]{Definition}
\newtheorem{question}[theorem]{Question}

\newcommand{\F}{\mathbb F}
\newcommand{\GP}{\mathrm{GP}}
\numberwithin{equation}{section}
\allowdisplaybreaks

\title{\bfseries Geometric progressions in distance sets in odd dimensions over prime fields}
\author{%
  Thang Pham\thanks{Institute of Mathematics and Interdisciplinary Sciences,
    Xidian University. \newline
\hspace*{0.45cm} Email: \texttt{thangphammath@xidian.edu.cn}.}%
  \and Dung The Tran \thanks{VNU University of Science, Hanoi, Vietnam. ~Email: {\tt tranthedung56@gmail.com}}
  \and Boqing Xue\thanks{Institute of Mathematical Sciences,
    ShanghaiTech University. Email: \texttt{xuebq@shanghaitech.edu.cn}.}%
}

\date{}

\begin{document}
\maketitle

\begin{abstract}
For $E\subset\F_p^d$, let
$\Delta(E)=\{Q_d(x-y):x,y\in E\}$, where
$Q_d(z)=z_1^2+\cdots+z_d^2$. In this paper, we introduce and study the problem of the existence of nondegenerate three-term geometric progressions in
the distance set $\Delta(E)$ when $d$ is odd. Let $M_{\GP,d}(p)$ be the largest
size of a set whose distance set contains no such progression. For every prime $p\geq7$ and every integer $m\geq0$, we have
\[
 p^mM_{\GP,1}(p)\leq M_{\GP,2m+1}(p)
 \leq4p^mM_{\GP,1}(p).
\]
For all sufficiently large primes $p$, we also prove
\[
 p^{m+\frac{7}{26}}\leq M_{\GP,2m+1}(p)
 \leq Cp^{m+\frac{16}{25}},
\]
where $C$ is an absolute positive constant.
\end{abstract}
\tableofcontents

\section{Introduction}\label{sec:main}

Throughout the paper, $p$ is an odd prime and $\F_p$ is the field
with $p$ elements. For $d\geq1$, let
\[
 Q_d(z)=z_1^2+\cdots+z_d^2,
 \qquad
 \Delta(E)=\{Q_d(x-y):x,y\in E\}
 \quad(E\subset\F_p^d).
\]
The finite-field Erd\H{o}s--Falconer problem asks how large a set
must be to determine a positive proportion of all possible
distances. Bourgain, Katz, and Tao~\cite{BKT2004} connected
planar distance estimates with the sum--product problem.
Iosevich and Rudnev~\cite{IosevichRudnev2007} developed a
Fourier-analytic approach and proved
that, for $d\geq2$, the condition
$|E|\geq Cp^{\frac{d+1}{2}}$ implies $|\Delta(E)|\geq cp$.
Hart, Iosevich, Koh, and Rudnev~\cite{HartIosevichKohRudnev2011}
showed that the exponent $\frac{d+1}{2}$ cannot in general be
lowered in odd dimensions. In even dimensions, it is conjectured that the sharp exponent is $d/2$. This was confirmed recently, up to a log factor, in \cite{PhamShenTranXue2026}. 

In this paper, we ask whether the structure of distance sets allows us to obtain
a prescribed configuration, such as an arithmetic or a geometric
progression, under weaker size assumptions on $E$ than those in the
distance estimates above. Previous work has considered square and
nonsquare distances~\cite{KohPhamPham2023}, and the values that can be represented as ratios of nonzero
distances~\cite{IosevichKohParshall2019,IosevichKohRakhmonov2024,pham2024}.

In particular, we take a first step in this direction by studying
three-term geometric progressions in distance sets. We start with the following definition.
\begin{definition}\label{def:GP}
Let $D\subset\F_p$. A \emph{nondegenerate three-term geometric
progression} in $D$ is a triple
\[
 x,\ xr,\ xr^2\in D,
 \qquad x,r\in\F_p^\times,\qquad r^2\ne1.
\]
\end{definition}
The three terms are nonzero and pairwise distinct. Equivalently,
they are distinct nonzero values $t_1,t_2,t_3$ satisfying
$t_2^2=t_1t_3$. Multiplication by a nonzero scalar preserves the
existence of such a progression. Unless otherwise stated, geometric progressions have three terms.

The main question is formulated formally as follows.

\begin{question}\label{question:GP}
How large must $E\subset\F_p^{d}$ be in order that
$\Delta(E)$ contain a nondegenerate geometric progression?
\end{question}

Positive-density results already give partial answers for the existence of such progressions. Szemer\'edi's theorem \cite{Szemeredi1975} states that, for each fixed $k\geq3$ and $\delta>0$, every subset of $\{1,\ldots,N\}$ of size at least $\delta N$ contains a $k$-term arithmetic progression when $N$ is sufficiently large. Applying this theorem to a short interval of exponents in the cyclic group $\F_p^\times$ gives the corresponding assertion for geometric progressions of every fixed length in sets of fixed positive density. For three terms, Roth's theorem
\cite{Roth1953} is sufficient. The quantitative bounds of Bloom and Sisask~\cite{BloomSisask2023} and Raghavan~\cite{Raghavan2026}
give more information: for every fixed $c,K>0$, the condition
\[
 |D|\geq\frac{cp}{(\log p)^K}
\]
forces such a progression in $D\subset\F_p$ when $p$ is
sufficiently large in terms of $c,K$.\footnote{For completeness, we give the short deduction in
Lemma~\ref{lem:positive-density-gp}.} Therefore, combining this conclusion with \cite{PhamShenTranXue2026}, the condition $|E|\gg p^{\frac{d}{2}}$ is sufficient when $d$ even. 
For $d=2m$, the exponent $m$ is sharp in the split cases, namely when
$m$ is even or $p\equiv1\pmod4$: there is a totally isotropic subspace
of size $p^m$ whose distance set is $\{0\}$.
In the nonsplit case, when $m$ is odd and $p\equiv3\pmod4$,
the optimal exponent remains open; see Section~\ref{sec:even-remarks}.

It is worth noting that the arithmetic and geometric properties of the distance set play an important role in the study of Question \ref{question:GP}. For arbitrary sets of field elements, a lower bound with a fixed
exponent below one does not force such a progression. More precisely, for each fixed $\varepsilon>0$, a size condition
$|D|\geq p^{1-\varepsilon}$ does not force a nondegenerate three-term geometric progression in $D\subset\F_p$. Indeed, applying Behrend's construction \cite{Behrend1946} to a short interval of exponents gives sets of size at least $p \, e^{-C\sqrt{\log p}}$ with no such progression; see Section~\ref{sec:remarks}.

We now turn to odd dimensions. For each $d\geq1$, define
\[
 M_{\GP,d}(p)=
 \max\{|E|:E\subset\F_p^d,\ \Delta(E)\text{ contains no
 nondegenerate geometric progression}\}.
\]
Our main result gives upper and lower bounds for this extremal quantity.

\begin{theorem}\label{thm:GP-odd}
There are absolute constants $C\geq1$ and $p_0\geq7$ such that,
for every prime $p\geq p_0$ and every integer $m\geq0$,
\begin{equation}\label{eq:odd-extremal-bounds}
 p^{m+\frac{7}{26}}
 \leq M_{\GP,2m+1}(p)
 \leq Cp^{m+\frac{16}{25}}.
\end{equation}
\end{theorem}

The constants are independent of the dimension. For $m\geq1$,
the sufficient exponent $m+\frac{16}{25}$ is strictly smaller
than the exponent $m+1=\frac{d+1}{2}$ for a positive proportion of distances
in odd dimensions.

For $A\subset\F_p$, we have $\Delta(A)=(A-A)^2$, where
$(A-A)^2=\{(a-b)^2:a,b\in A\}$.

The main ingredient in the proof is the following one-dimensional result.

\begin{proposition}\label{prop:one-dimensional}
The following statements hold.
\begin{enumerate}
\renewcommand{\labelenumi}{\textup{(\roman{enumi})}}
\item There is an absolute constant $C$ such that, for every prime
      $p\geq7$ and every $A\subset\F_p$, the condition
      \[
       |A|\geq Cp^{\frac{16}{25}}
      \]
      implies that $(A-A)^2$ contains a nondegenerate geometric
      progression.
\item For every sufficiently large prime $p$, there exists a set
      $A\subset\F_p$ such that
      \[
       |A|\geq p^{\frac{7}{26}}
      \]
      and $\kappa(A-A)^2$ contains no nondegenerate geometric
      progression for every $\kappa\in\F_p^\times$.
\end{enumerate}
\end{proposition}

Geometric progressions in sumsets provide another relevant comparison.
Ahmadi and Shparlinski~\cite[Theorem~1]{AhmadiShparlinski2007} proved that
$A+B$ contains three distinct nonzero terms in geometric
progression whenever $A,B\subset\F_p$ and
$|A||B|\geq20p^{\frac43}$.
Taking $B=-A$ gives such a progression in $A-A$ when
$|A|\geq\sqrt{20}\,p^{\frac23}$.
Shparlinski~\cite{Shparlinski2020} strengthened the vector-space
result and obtained progressions with a separate ratio in each
coordinate.

Proposition~\ref{prop:one-dimensional}\textup{(i)} gives the following improvement for difference sets.

\begin{corollary}\label{cor:difference-set-gp}
There is an absolute constant $C$ such that, for every prime
$p\geq7$ and every $A\subset\F_p$ with
$|A|\geq Cp^{\frac{16}{25}}$, there exist $x,r\in\F_p^\times$ such that
\[
 x,\ xr,\ xr^2\in A-A,\qquad r^4\ne1.
\]
\end{corollary}

\paragraph{Sketch of main ideas of the proofs.}
We first reduce the problem in odd dimensions to dimension one, using a one-dimensional version of the quotient argument of
Pham, Shen, and Xue
\cite{PhamShenXue2026}.
We choose an $m$-dimensional totally isotropic subspace $W$. More precisely, averaging over its orthogonal-group orbit controls the number
of pairs of points of $E$ in common affine cosets of $W$.
In a suitable affine coset of $W^\perp$, projection to the
one-dimensional space $W^\perp/W$ then produces a large set
$A\subset\F_p$ with
$\kappa(A-A)^2\subseteq\Delta(E)$ for some $\kappa\ne0$.
For the lower bound, replacing each point of $A$ by a
translate of $W$ multiplies the cardinality by $p^m$ and
changes the distance set only by a nonzero scalar.

To prove Proposition~\ref{prop:one-dimensional}\textup{(i)},
suppose that $(A-A)^2$
contains no nondegenerate geometric progression. For
$t\ne0$ with $t^4\ne1$, the map
\[
 (a,b,c)\longmapsto(b-ta,c-tb)
\]
is injective on $A^3$. This gives $|A|^3\leq p^2$ and bounds
the sizes of intersections of translates and dilates of $A$.
Combining this injectivity with Pl\"unnecke's inequality gives
$|A-A|\geq|A|^{\frac54}$.
We apply these estimates to the variance of the number of
representations as a product of two nonzero differences.
At every nonzero squared difference this number is small,
whereas its mean is large when $|A|$ is large. Comparing the
resulting lower and upper bounds for the variance gives
$|A|\ll p^{\frac{16}{25}}$.

For part~\textup{(ii)}, we select vectors in $\{0,1,2,3\}^L$ with
prescribed digit counts and project them by a random linear map to $\F_p$.
We delete points to remove relations that would produce
nondegenerate geometric progressions among squared differences.
The choice of four digits ensures that polynomial identities
give only degenerate progressions. Suitable digit weights
bound the number of deleted points and leave a set of size
at least $p^{\frac{7}{26}}$.
The cylinder construction then gives the examples
in every odd dimension.

{\bf Structure.} Section~\ref{sec:reduction} contains the reduction to dimension one
and deduces Theorem~\ref{thm:GP-odd} from
Proposition~\ref{prop:one-dimensional}.
Sections~\ref{sec:one-dimensional-upper} and~\ref{sec:odd-examples}
prove parts~\textup{(i)} and~\textup{(ii)} of that proposition,
respectively. Section~\ref{sec:remarks} discusses even dimensions
and the remaining extremal questions.

We write $X\ll Y$ if $X\leq CY$ for an absolute constant $C>0$,
$X\gg Y$ if $Y\ll X$, and $X\asymp Y$ if both inequalities
hold. Constants may vary from one occurrence to the next.
We use $\log$ for the natural logarithm.

\section{Reduction: from \texorpdfstring{$2m+1$ to $1$}{2m+1 to 1}}\label{sec:reduction}

We compare the extremal problem in odd dimensions with the
one-dimensional problem. The comparison does not depend on the
exponents in Theorem~\ref{thm:GP-odd}.

\begin{proposition}\label{prop:odd-dimension-comparison}
For every prime $p\geq7$ and every integer $m\geq0$,
\begin{equation}\label{eq:finite-comparison}
 p^mM_{\GP,1}(p)\leq M_{\GP,2m+1}(p)
 \leq4p^mM_{\GP,1}(p).
\end{equation}
\end{proposition}

The constants are independent of the dimension. Thus, every
one-dimensional upper bound or construction transfers to all odd
dimensions, with the factor $p^m$.
For the upper bound, we extract a large set in an isotropic quotient.
For the lower bound, we replace each point of a one-dimensional set
by a translate of a totally isotropic subspace.
We first record the quadratic-form and Fourier identities used in
these two arguments.

\subsection{Quadratic forms and isotropic quotients}\label{sec:quadratic-quotients}

Let $Q$ be a nondegenerate quadratic form on $\F_p^d$.
Its associated bilinear form is
\[
 B_Q(u,v)=\frac{Q(u+v)-Q(u)-Q(v)}2.
\]
For a subspace $W\subset\F_p^d$, we write
\[
 W^\perp=\{u\in\F_p^d:B_Q(u,w)=0\text{ for every }w\in W\}.
\]
For the standard form $Q_d$, we have $B_{Q_d}(u,v)=u\cdot v$.
All orthogonal complements below are taken with respect to the form under consideration.

The subspace $W$ is \emph{totally isotropic} if $Q(w)=0$ for every
$w\in W$. Since $p$ is odd, this implies $B_Q(w,w')=0$ for all
$w,w'\in W$. Hence, $W\subset W^\perp$. If $u\in W^\perp$ and
$w\in W$, then
\[
 Q(u+w)=Q(u)+2B_Q(u,w)+Q(w)=Q(u).
\]
Thus, the value of $Q(u)$ depends only on the coset $u+W$, so $Q$
defines a quadratic form on $W^\perp/W$.
In particular, passing from an affine coset of $W^\perp$ to
this quotient preserves the quadratic values of all differences.

The induced form is nondegenerate.
Indeed, if $u\in W^\perp$ satisfies $B_Q(u,v)=0$ for every $v\in W^\perp$,
then $u\in(W^\perp)^\perp=W$, so its class in the quotient is zero.
Since $\dim W^\perp=d-\dim W$, we also have
\begin{equation}\label{eq:quotient-dimension}
 \dim(W^\perp/W)=d-2\dim W.
\end{equation}

Witt decomposition shows that every nondegenerate quadratic
form on $\F_p^{2m+1}$ has a totally isotropic subspace of
dimension $m$; see \cite[Lemma~2.1]{IosevichKohRakhmonov2024}.
For such a subspace $W$, the quotient $W^\perp/W$ is
one-dimensional and nondegenerate. Therefore, the induced quadratic form is $\kappa t^2$ with $\kappa\in\F_p^\times$.

\subsection{Extraction of a one-dimensional set}

The next lemma extracts a set $A\subset\F_p$ whose squared
differences, after multiplication by a nonzero scalar, belong to
$\Delta(E)$.

\begin{lemma}\label{lem:odd-projection}
Let $p\geq7$, let $m\geq1$, and let
$\varnothing\ne E\subset\F_p^{2m+1}$. There are $A\subset\F_p$
and $\kappa\in\F_p^\times$ such that
\begin{equation}\label{eq:projection-size}
 |A|\geq\frac{|E|}{p^m(2+\frac{|E|}{p^{m+1}})},
 \qquad \kappa(A-A)^2\subseteq\Delta(E).
\end{equation}
\end{lemma}

When $|E|\leq p^{m+1}$, the denominator in parentheses in
\eqref{eq:projection-size} is at most $3$. We retain the full
bound so that Proposition~\ref{prop:odd-dimension-comparison}
can be proved without an initial restriction on $|E|$.

\begin{proof}
We seek a totally isotropic $m$-dimensional subspace $W$ and an
affine coset of $W^\perp$ in which many cosets of $W$ meet $E$.
After translation, their images in $W^\perp/W$ will form the
required set $A$.
To find $W$, we estimate the number of pairs of points of $E$
lying in the same coset of $W$.

We first bound the number of pairs at distance zero.
Write $S_0=\{x\in\F_p^{2m+1}:Q_{2m+1}(x)=0\}$. 
We have $|S_0|=p^{2m}$ and
\begin{equation}\label{eq:zero-incidence}
 |\{(x,y)\in E^2: Q_{2m+1}(x-y)=0\}|
 \leq \frac{|E|^2}{p}+p^m|E|,
\end{equation}
as follows from the sphere Fourier formula in
\cite[Section~4.1]{IosevichKohLewko2020}.
We sketch the argument to make the constant explicit.
For $f:\F_p^{2m+1}\to\mathbb C$, write

\[
 \widehat f(\xi)=p^{-(2m+1)}\sum_{x\in\F_p^{2m+1}}
 f(x)e^{-2\pi i x\cdot\xi/p}.
\]
The sphere formula has a remaining character sum that is either zero
or has absolute value $\sqrt p$. Together with Plancherel's identity, this give
\[
 \begin{aligned}
 \widehat{\mathbf1_{S_0}}(0)&=p^{-1},\qquad
 |\widehat{\mathbf1_{S_0}}(\xi)|\leq p^{-m-1}\quad(\xi\ne0),\\
 \sum_{\xi\in\F_p^{2m+1}}|\widehat{\mathbf1_E}(\xi)|^2
 &=p^{-(2m+1)}|E|.
 \end{aligned}
\]
Fourier inversion therefore yields
\[
 \begin{aligned}
 \sum_{x,y\in E}\mathbf1_{S_0}(x-y)
 &=p^{4m+2}\sum_{\xi\in\F_p^{2m+1}}
   \widehat{\mathbf1_{S_0}}(\xi)|\widehat{\mathbf1_E}(\xi)|^2\\
 &\leq\frac{|E|^2}{p}
   +p^{4m+2}p^{-m-1}p^{-(2m+1)}|E|\\
 &=\frac{|E|^2}{p}+p^m|E|,
 \end{aligned}
\]
which proves \eqref{eq:zero-incidence}.

We now average over isotropic subspaces. Choose one such subspace
of dimension $m$, and choose $W$ uniformly from its orbit under
the orthogonal group. To compute the probability that a fixed
vector lies in $W$, we use transitivity on the nonzero isotropic
vectors. 
For $Q_{2m+1}(u)\ne0$, reflection across $u^\perp$ is the map
\[
 x\longmapsto x-\frac{2x\cdot u}{Q_{2m+1}(u)}u.
\]
For two nonzero isotropic vectors
$v,v'$,
reflection across the hyperplane orthogonal to $v-v'$ sends $v$ to $v'$ if $v\cdot v'\ne0$.
Indeed, $Q_{2m+1}(v-v')=-2v\cdot v'$ and
$v\cdot(v-v')=-v\cdot v'$.
If $v\cdot v'=0$, choose $u$ with
$Q_{2m+1}(u)(v\cdot u)(v'\cdot u)\ne0$ and first reflect $v$ across the hyperplane orthogonal to $u$.
Its inner product with $v'$ is then nonzero, since
\[
 \left(v-\frac{2v\cdot u}{Q_{2m+1}(u)}u\right)\cdot v'
 =-\frac{2(v\cdot u)(v'\cdot u)}{Q_{2m+1}(u)}\ne0.
\]
A second reflection therefore sends this image of $v$ to $v'$.
Such a $u$ exists
because the  nonzero polynomial
$u\mapsto Q_{2m+1}(u)(v\cdot u)(v'\cdot u)$ has degree four, less
than $p$, and so cannot vanish at every point of $\F_p^{2m+1}$.
It follows by double counting that, for every nonzero isotropic
vector $v$,
\begin{equation}\label{eq:orbit-probability}
 \Pr(v\in W)=\frac{p^m-1}{p^{2m}-1}
            =\frac1{p^m+1}<p^{-m}.
\end{equation}

The sum $\sum\limits_{x+W}|E\cap(x+W)|^2$ counts ordered pairs of
points whose difference lies in $W$. By
\eqref{eq:zero-incidence} and \eqref{eq:orbit-probability}, its
expectation is at most
\[
 |E|+p^{-m}\left(\frac{|E|^2}{p}+p^m|E|\right)
 =2|E|+\frac{|E|^2}{p^{m+1}}.
\]
Fix $W$ for which this bound holds. The Cauchy--Schwarz inequality
then shows that at least
\[
 \frac{|E|}{2+\frac{|E|}{p^{m+1}}}
\]
affine cosets of $W$ meet $E$. Each lies in one affine coset of
$W^\perp$. Since there are $p^m$ cosets of $W^\perp$, some affine subspace $z+W^\perp$ contains at least
\[
 \frac{|E|}{p^m(2+\frac{|E|}{p^{m+1}})}
\]
cosets of $W$ that meet $E$.

Translate $E\cap(z+W^\perp)$ by $-z$ and take its image in
$W^\perp/W$.
By \eqref{eq:quotient-dimension}, the quotient is one-dimensional
and nondegenerate, so its form can be written as $\kappa t^2$ with
$\kappa \ne0$. Identify the image with a set $A\subset\F_p$, and
choose one representative $x_a\in E\cap(z+W^\perp)$ for each
$a\in A$. The quotient form gives
\[
\kappa(a-b)^2= Q_{2m+1}(x_a-x_b)\qquad(a,b\in A).
\]
This proves both assertions and completes the proof.
\end{proof}

\subsection{The cylinder construction}

The reverse construction replaces each element of $A$ by a
translate of an $m$-dimensional totally isotropic subspace.
This multiplies the cardinality by $p^m$, while the distances
depend only on the corresponding elements of $A$. The following
lemma is the coordinate-free form of the construction in
\cite[Lemma~2.1(ii) and Section~4.2]{IosevichKohRakhmonov2024}.

\begin{lemma}
\label{lem:Witt-cylinder}
Let $m\geq0$ be an integer and let $d=2m+1$.
There exists a totally isotropic $m$-dimensional subspace
$W\subset\F_p^d$ and a vector $e\in W^\perp\setminus W$
such that $\kappa:=Q_d(e)\ne0$ and,
for every $A\subset\F_p$, the set
\[
 E_A:=\{w+ae:w\in W,\ a\in A\}
\]
satisfies
\begin{equation}
 |E_A|=p^m|A|,
 \qquad
 \Delta(E_A)=\kappa(A-A)^2.
 \label{eq:Witt-cylinder}
\end{equation}
\end{lemma}

\begin{proof}
If $m=0$, then the conclusion follows by taking $W=\{0\}$ and $e=1$, so that
$\kappa=Q_1(e)=1$. Therefore, it remains to prove the result in the case
$m\geq1$.

Choose a totally isotropic subspace $W\subset\F_p^d$ of dimension $m$, as in
Section~\ref{sec:quadratic-quotients}. The quadratic form induced by $Q_d$ on
$W^\perp/W$ is nondegenerate and one-dimensional. Hence, for any
$e\in W^\perp\setminus W$,
\[
\kappa:=Q_d(e)\ne0.
\]

Since $e\notin W$, the representation $w+ae$, with $w\in W$ and $a\in A$, is unique. Therefore,
\[
|E_A|=|W||A|=p^m|A|.
\]
Moreover, for $w,w'\in W$ and $a,b\in A$, the total isotropy of $W$ and the fact that $e\in W^\perp$ give
\[
Q_d(w-w')=0,
\qquad
(w-w')\cdot e=0.
\]
Consequently,
\[
\begin{aligned}
Q_d\bigl((w+ae)-(w'+be)\bigr)
&=Q_d(w-w')+2(a-b)(w-w')\cdot e\\
&\quad +(a-b)^2Q_d(e)\\
&=\kappa(a-b)^2.
\end{aligned}
\]
Thus, $\Delta(E_A)\subseteq\kappa(A-A)^2$. 
The reverse inclusion follows by taking $w=w'=0$. Hence,
$\Delta(E_A)=\kappa(A-A)^2$. This completes the proof.
\end{proof}

\subsection{Proof of Proposition~\ref{prop:odd-dimension-comparison}, Theorem~\ref{thm:GP-odd}, and Corollary \ref{cor:difference-set-gp}}

We now combine the two lemmas. The estimate
$M_{\GP,1}(p)\leq(p-1)/2$ allows us to absorb the dependence
on $|E|$ in the denominator of \eqref{eq:projection-size}.

\begin{proof}[Proof of Proposition~\ref{prop:odd-dimension-comparison}]
If $m=0$, then \eqref{eq:finite-comparison} follows directly from the definition of $M_{\GP,1}(p)$. We henceforth assume that $m\geq1$.
We first note that
\[
 M_{\GP,1}(p)\leq\frac{p-1}{2}.
\]
Indeed, if $A\subset\F_p$ has more than $\frac{p}{2}$ elements,
then $A$ and $u+A$ intersect for every $u\in\F_p$, so
$A-A=\F_p$. Hence, $(A-A)^2$ contains all nonzero squares.
The set of nonzero squares is a cyclic multiplicative subgroup of order $\frac{p-1}{2}\geq3$. Therefore, it contains a nondegenerate
geometric progression.

Let $E\subset\F_p^{2m+1}$ be nonempty and suppose that its
distance set contains no nondegenerate geometric progression.
By Lemma~\ref{lem:odd-projection} and the invariance under nonzero
scaling noted after Definition~\ref{def:GP},
\[
 M_{\GP,1}(p)
 \geq\frac{|E|}{p^m\left(2+\frac{|E|}{p^{m+1}}\right)}.
\]
Multiplying through by the denominator and collecting the terms
containing $|E|$, we obtain
\[
 |E|\left(1-\frac{M_{\GP,1}(p)}{p}\right)
 \leq2p^mM_{\GP,1}(p).
\]
Since $M_{\GP,1}(p)\leq\frac{p-1}{2}$, the factor on the left is greater than
$\frac{1}{2}$. Hence, $|E|\leq4p^mM_{\GP,1}(p)$.  This proves the upper bound in \eqref{eq:finite-comparison}.

For the lower bound, choose $A\subset\F_p$ with
$|A|=M_{\GP,1}(p)$ such that $(A-A)^2$ contains no nondegenerate
geometric progression. Lemma~\ref{lem:Witt-cylinder} constructs a set $E_A\subset\F_p^{2m+1}$ satisfying
$|E_A|=p^mM_{\GP,1}(p)$ and
$\Delta(E_A)=\kappa(A-A)^2$ for some $\kappa\ne0$.
Since multiplication by $\kappa$ preserves nondegeneracy,
$\Delta(E_A)$ also contains no nondegenerate geometric progression. Consequently,
$M_{\GP,2m+1}(p)\geq p^mM_{\GP,1}(p)$.
This proves \eqref{eq:finite-comparison} and completes the proof.
\end{proof}

The main theorem follows by applying the comparison to the two
one-dimensional bounds.

\begin{proof}[Proof of Theorem~\ref{thm:GP-odd}]
Proposition~\ref{prop:one-dimensional}\textup{(i)} gives an
absolute constant $C_1$ such that
$M_{\GP,1}(p)\leq C_1p^{\frac{16}{25}}$ for every prime $p\geq7$.
Part~\textup{(ii)} gives an absolute $p_0\geq7$ such that
$M_{\GP,1}(p)\geq p^{\frac{7}{26}}$ for every prime $p\geq p_0$.
Hence, Proposition~\ref{prop:odd-dimension-comparison} yields
\[
 p^{m+\frac{7}{26}}
 \leq p^mM_{\GP,1}(p)
 \leq M_{\GP,2m+1}(p)
 \leq4p^mM_{\GP,1}(p)
 \leq4C_1p^{m+\frac{16}{25}}.
\]
Taking $C=\max\{1,4C_1\}$ proves the theorem. Both $C$ and
$p_0$ are independent of $m$. This completes the proof.
\end{proof}

\begin{proof}[Proof of Corollary \ref{cor:difference-set-gp}]
By Proposition~\ref{prop:one-dimensional}\textup{(i)}, there exist
$u,v,w\in(A-A)\setminus\{0\}$ such that $u^2,v^2,w^2$ are
pairwise distinct and form a geometric progression. Thus,
\[
v^4=u^2w^2,
\]
and hence
\[
(v^2-uw)(v^2+uw)=0.
\]
Therefore, $v^2=uw$ or $v^2=-uw$. Since $A-A$ is symmetric,
replacing $w$ by $-w$ if necessary, we may assume that
$v^2=uw$.

Set $x=u$ and $r=v/u$. Then $x,r\in\F_p^\times$ and
\[
x=u,\qquad xr=v,\qquad xr^2=\frac{v^2}{u}=w.
\]
Thus, $x,xr,xr^2\in A-A$. Moreover,
\[
r^4=\frac{v^4}{u^4}=\frac{w^2}{u^2}\ne1,
\]
since $u^2\ne w^2$. This completes the proof.
\end{proof}

We prove the two parts of Proposition~\ref{prop:one-dimensional}
in the next two sections.

\section{Proof of Proposition~\ref{prop:one-dimensional}\textup{(i)}}\label{sec:one-dimensional-upper}
The proof compares two bounds for the variance of the number of
representations of a nonzero field element as a product of two
differences from $A$. The absence of a progression makes this
number small at every nonzero squared difference. We first show
that there are many such squared differences and establish the
injectivity estimate needed for the upper bound on the variance.

\subsection{Consequences of the absence of a progression}

The three conclusions of the next lemma have different roles.
The bound for $|A|$ controls the exceptional slopes in the variance
estimate, the bound for $|A-A|$ supplies many squared differences,
and the bound for the number of solutions to $b-ta=s$ controls
the remaining slopes. All three follow from the injectivity of
$(a,b,c)\mapsto(b-ta,c-tb)$ when $t\ne0$ and $t^4\ne1$.

\begin{lemma}\label{lem:hereditary}
Let $p\geq7$ be a prime and let $\varnothing\ne A\subset\F_p$.
If $(A-A)^2$
contains no nondegenerate geometric progression, then
\begin{equation}\label{eq:new-difference-growth}
 |A|^3\leq p^2,
 \qquad |A-A|\geq|A|^{\frac54}.
\end{equation}
Moreover, for $t\ne0$ with $t^4\ne1$ and every $s\in\F_p$,
\begin{equation}\label{eq:affine-fiber-cap}
 |\{(a,b)\in A^2:b-ta=s\}|\leq\frac p{|A|}.
\end{equation}
\end{lemma}
\begin{proof}
Fix $t\in\F_p^\times$ with $t^4\ne1$ and consider the map
\[
 \Phi_t:A^3\longrightarrow\F_p^2,
 \qquad \Phi_t(a,b,c)=(b-ta,c-tb).
\]
If two triples $(a,b,c)$ and $(a',b',c')$ have the same image, then
\[
 b-b'=t(a-a'),\qquad c-c'=t(b-b').
\]
If the triples are distinct, then $a-a'\ne0$, and these three
differences are nonzero and have common ratio $t$. Their squares
would form a nondegenerate geometric progression in $(A-A)^2$,
since $t^4\ne1$. Thus, $\Phi_t$ is injective.
Since $2^4\ne1$ for $p\geq7$, taking $t=2$ yields
$|A|^3\leq p^2$.

To obtain the lower bound for $|A-A|$, apply Pl\"unnecke's
inequality \cite[Theorem~1.1]{Petridis2012} to $A$ and $-A$.
This gives a a nonempty set $X\subset A$ such that
\[
 |X-A-A|\leq
 \left(\frac{|A-A|}{|A|}\right)^2|X|.
\]
For every $(a,b,c)\in A\times X\times X$, we have
\[
 b-2a=b-a-a\in X-A-A,
 \qquad
 c-2b=c-b-b\in X-A-A.
\]
Since $\Phi_2$ is injective,
\[
 |A||X|^2\leq|X-A-A|^2
 \leq\left(\frac{|A-A|}{|A|}\right)^4|X|^2.
\]
Since $X\ne\varnothing$, it follows that
$|A-A|\geq|A|^{\frac{5}{4}}$.

Finally, for each $s\in\F_p$, the map $\Phi_t$ sends
$\{(a,b,c)\in A^3:b-ta=s\}$ injectively into
$\{s\}\times\F_p$. Comparing cardinalities gives
\[
 |A|\,|\{(a,b)\in A^2:b-ta=s\}|\leq p,
\]
which proves \eqref{eq:affine-fiber-cap}. This completes the proof.
\end{proof}

\subsection{The variance argument}

For a nonempty set $A\subset\F_p$, let
$r_{A-A}(h)=|\{(a,b)\in A^2:a-b=h\}|$ for $h\in\F_p$.
Thus, $r_{A-A}(h)>0$ if and only if $h\in A-A$.
To count products of nonzero differences with their multiplicities, define
\[
 F(z):=\sum_{\substack{x,y\in\F_p^\times}\\xy=z} r_{A-A}(x)r_{A-A}(y)
 \quad(z\in\F_p^\times).
\]
Its mean is
\[
 \overline F:=\frac{1}{p-1}\sum_{z\in\F_p^\times}F(z)
 =\frac{|A|^2(|A|-1)^2}{p-1}.
\]
The formula for the mean follows from
$\sum_{h\ne0} r_{A-A}(h)=|A|(|A|-1)$.

Let $h\in(A-A)\setminus\{0\}$. If $x,y\in\F_p^\times$
satisfy $xy=h^2$ and $r_{A-A}(x)r_{A-A}(y)>0$, then
$x,y\in A-A$. Consequently, $x^2,h^2,y^2\in(A-A)^2$ and
\[
 \frac{h^2}{x^2}=\frac{y}{x}=\frac{y^2}{h^2}.
\]
Thus, these three values form a geometric progression, possibly degenerate.
When no nondegenerate progression is present, only a bounded
number of pairs $(x,y)$ are possible. Thus, $F(h^2)$ is small.
For $|A|>6\sqrt p$, these values lie below half the mean and give
the lower bound in the following lemma. The upper bound will
follow from Lemma~\ref{lem:hereditary}.

\begin{lemma}\label{lem:product-variance-bounds}
Let $p\geq7$ be a prime and let $\varnothing\ne A\subset\F_p$.
Suppose that $(A-A)^2$ contains no nondegenerate geometric
progression. Then
\[
 \sum_{z\in\F_p^\times}(F(z)-\overline F)^2
 \leq5p^2|A|^3.
\]
If $|A|>6\sqrt p$, then also
\[
 \sum_{z\in\F_p^\times}(F(z)-\overline F)^2
 \gg\frac{(|A-A|-1)|A|^8}{p^2},
\]
where the implied constant is absolute.
\end{lemma}

\begin{proof}[Proof of Lemma~\ref{lem:product-variance-bounds}]
We first bound the variance from above by counting pairs on lines
of a fixed slope. We then obtain the lower bound from the small
values of $F$ at nonzero squared differences.

The equation $xy=uv\ne0$ can be written as $u=tx$ and $y=tv$
for some $t\ne0$. This relates the second moment of $F$ to the
numbers of solutions of $b-b'=t(a-a')$. We count these solutions
through the fibers of $(a,b)\mapsto b-ta$. For $s,t\in\F_p$, let
\[
 i_t(s)=|\{(a,b)\in A^2:b-ta=s\}|,
 \qquad V(t)=\sum_{s\in\F_p}\left(i_t(s)-\frac{|A|^2}{p}\right)^2.
\]
Since $\sum\limits_{s\in\F_p} i_t(s)=|A|^2$, expanding the square in $V(t)$ reduces the calculation to counting solutions of
$b-b'=t(a-a')$. When $a\ne a'$, each choice of $a,a',b,b'$
gives one value of $t$; when $a=a'$, a solution requires $b=b'$
and then all $p$ values of $t$ are allowed. Hence,
\[
 \sum_{t\in\F_p}\sum_{s\in\F_p} i_t(s)^2=(|A|^2-|A|)|A|^2+p|A|^2.
\]
Subtracting $|A|^4$ gives
\begin{equation}\label{eq:new-total-variance}
 V(t)=\sum_{s\in\F_p} i_t(s)^2-\frac{|A|^4}{p}\geq0,
 \qquad \sum_{t\in\F_p}V(t)=p|A|^2-|A|^3\leq p|A|^2.
\end{equation}
For $t\ne0$ with $t^4\ne1$, \eqref{eq:affine-fiber-cap} implies
$V(t)\leq\sum\limits_{s\in\F_p} i_t(s)^2\leq p|A|$.
For the at most four values of $t$ satisfying $t^4=1$, we use
$V(t)\leq\sum\limits_{s\in\F_p} i_t(s)^2\leq|A|^3$. 
For $t\ne0$ with $t^4\ne1$, we have $V(t)^2\leq p|A|V(t)$.
The values satisfying $t^4=1$ contribute at most $4|A|^6$.
Thus, \eqref{eq:new-total-variance} and the bound $|A|^3\leq p^2$
from Lemma~\ref{lem:hereditary} yield
\begin{equation}\label{eq:new-variance-upper}
 \sum_{t\in\F_p^\times}V(t)^2
 \leq p^2|A|^3+4|A|^6
 \leq 5p^2|A|^3.
\end{equation}
We now transfer this estimate to $F$. Grouping the weighted solutions of $xy=uv\ne0$ by $t=\frac{u}{x}=\frac{y}{v}$ gives
\begin{equation}\label{eq:product-ratio-second-moment}
 \sum_{z\in\F_p^\times}F(z)^2
 =\sum_{t\in\F_p^\times}
   \left(\sum_{h\in\F_p^\times} r_{A-A}(h)r_{A-A}(th)\right)^2.
\end{equation}
Moreover, since multiplication by any $h\ne0$ permutes $\F_p^\times$,
\[
 \sum_{t\in\F_p^\times}\sum_{h\in\F_p^\times}r_{A-A}(h)r_{A-A}(th)
 =\left(\sum_{h\in\F_p^\times}r_{A-A}(h)\right)^2
 =|A|^2(|A|-1)^2
 =\sum_{z\in\F_p^\times}F(z).
\]
Thus, $F$ and the function
$t\mapsto\sum\limits_{h\in\F_p^\times}r_{A-A}(h)r_{A-A}(th)$
have the same mean $\overline F$ and second moment.

For $t\in\F_p^\times$, the sum $\sum\limits_{s\in\F_p}i_t(s)^2$
counts solutions of $b-b'=t(a-a')$ in $A^4$.
Removing the $|A|^2$ solutions with $a=a'$ gives
\[
 \sum_{h\in\F_p^\times}r_{A-A}(h)r_{A-A}(th)
 =\sum_{s\in\F_p}i_t(s)^2-|A|^2
 =V(t)+\frac{|A|^4}{p}-|A|^2.
\]
Combining this identity with
\eqref{eq:product-ratio-second-moment} and the equality of the means,
we obtain
\begin{align}
 \sum_{z\in\F_p^\times}(F(z)-\overline F)^2
 &=\sum_{t\in\F_p^\times}
   \left(V(t)-\frac{1}{p-1}
   \sum_{u\in\F_p^\times}V(u)\right)^2\nonumber\\
 &\leq\sum_{t\in\F_p^\times}V(t)^2
 \leq 5p^2|A|^3,
 \label{eq:new-product-variance}
\end{align}
where the last inequality follows from \eqref{eq:new-variance-upper}.

We now prove the lower bound. For $h\in(A-A)\setminus\{0\}$, each nonzero summand in $F(h^2)$
corresponds to nonzero $x,y\in A-A$ with $xy=h^2$.
The terms $x^2,h^2,y^2$ form a geometric progression with ratio
$\frac{h^2}{x^2}$. This progression must be degenerate, so
$\left(\frac{h^2}{x^2}\right)^2=1$, or equivalently $x^4=h^4$.
There are at most four choices for $x$, each determining $y$.
Since $r_{A-A}\leq|A|$, we obtain
\begin{equation}\label{eq:new-product-low}
 F(h^2)\leq4|A|^2.
\end{equation}
Assume now that $|A|>6\sqrt p$. Then
\[
 \overline F
 =\frac{|A|^2(|A|-1)^2}{p-1}
 \geq\frac{|A|^4}{4p}.
\]
Since $|A|>6\sqrt p$, we also have $\overline F>8|A|^2$.

Since $A-A$ is symmetric and $p$ is odd, the squaring map is
two-to-one on $(A-A)\setminus\{0\}$, so the set
$D=(A-A)^2\setminus\{0\}$ has cardinality
$|D|=\frac{|A-A|-1}{2}$.
Moreover, \eqref{eq:new-product-low} and the bound
$\overline F\geq8|A|^2$ give
$F(z)\leq4|A|^2\leq\frac{\overline F}{2}$ for every $z\in D$. Therefore,
\begin{equation}\label{eq:variance-lower-explicit}
\begin{aligned}
 \sum_{z\in\F_p^\times}(F(z)-\overline F)^2
 &\geq\sum_{z\in D}(F(z)-\overline F)^2\\
 &\geq\frac{|D|}{4}\,\overline F^{\,2}
  =\frac{|A-A|-1}{8}\,\overline F^{\,2}\\
 &\gg\frac{(|A-A|-1)|A|^8}{p^2}.
\end{aligned}
\end{equation}

This completes the proof.
\end{proof}

We can now prove the upper bound. The next calculation explains
the exponent $\frac{16}{25}$.

\begin{proof}[Proof of Proposition~\ref{prop:one-dimensional}\textup{(i)}]
Assume that $(A-A)^2$ contains no nondegenerate geometric
progression. If $|A|\leq6\sqrt p$, then the required bound
$|A|\ll p^{\frac{16}{25}}$ is immediate. Hence, we may assume
$|A|>6\sqrt p$.
Comparing the two estimates in Lemma~\ref{lem:product-variance-bounds}, we obtain
\[
 (|A-A|-1)|A|^5\ll p^4.
\]
Since \eqref{eq:new-difference-growth} and $|A|>6\sqrt p$ imply
$|A-A|-1\gg|A|^{\frac{5}{4}}$, it follows that
\[
 |A|^{\frac{25}{4}}\ll p^4,
 \qquad\text{so}\qquad
 |A|\ll p^{\frac{16}{25}}.
\]
Choosing the constant in Proposition~\ref{prop:one-dimensional}\textup{(i)}
sufficiently large proves the assertion. This completes the proof.
\end{proof}

\section{Proof of Proposition~\ref{prop:one-dimensional}\textup{(ii)}}\label{sec:odd-examples}

In this section, we construct a large family of integer vectors, project it to
$\F_p$, and delete points to remove the relations that would
produce a nondegenerate geometric progression.

\subsection{Random projection and deletion}

To construct the set in Proposition~\ref{prop:one-dimensional}\textup{(ii)},
we must control relations among all differences of $A$.
The following criterion allows us to work before squaring
those differences.

For $\kappa\in\F_p^\times$, the set $\kappa(A-A)^2$ contains a
nondegenerate geometric progression if and only if there exist
$u,v,w\in(A-A)\setminus\{0\}$ with
\begin{equation}\label{eq:GP-difference-criterion}
 v^2=uw,\qquad u^4\ne v^4.
\end{equation}
Indeed, if $v^2=uw$, then $v^4=u^2w^2$, so
$\kappa u^2,\kappa v^2,\kappa w^2$ form a geometric progression
with common ratio $v^2/u^2$. All three terms are nonzero, and
$u^4\ne v^4$ says exactly that the square of this ratio is not one.
Conversely, a nondegenerate progression $\kappa u^2,\kappa v^2,\kappa w^2$ satisfies $v^4=u^2w^2$ and $u^4\ne v^4$. Since $p$ is odd, $v^2=uw$ or $v^2=-uw$. In the second case, replace $w$ with $-w$; this is allowed because $A-A$ is symmetric, and it leaves $w^2$ unchanged. Thus, \eqref{eq:GP-difference-criterion} holds.

We begin with a family of integer vectors in $\{0,1,2,3\}^L$
and project it linearly to $\F_p$. A \textit{collision} in the projection or an unwanted \textit{quadratic relation} can be removed by deleting points so that one of its differences no longer occurs.

\begin{lemma}\label{lem:digit-projection}
Let $L\geq1$ and $\mathcal C\subseteq\{0,1,2,3\}^L$. With differences and intersections taken in $\mathbb Z^L$, let
\[
 \mathcal D:=\mathcal C-\mathcal C,\qquad
r_{\mathcal C-\mathcal C}(\mathbf u):=|\mathcal C\cap(\mathcal C+\mathbf u)|,
 \qquad
 J(\mathcal C):=\sum_{\mathbf u\in\mathcal D}r_{\mathcal C-\mathcal C}(\mathbf u)^{\frac13}.
\]
For every prime $p>18$, there exists $A\subset\F_p$ such that
$(A-A)^2$ contains no nondegenerate geometric progression and
\begin{equation}\label{eq:deletion-projection}
 |A|\geq|\mathcal C|-\frac{|\mathcal C|^2}{p}-\frac{2J(\mathcal C)^3}{p}.
\end{equation}
\end{lemma}

\begin{proof}
We seek a projection that becomes injective after a controlled
number of deletions and creates no nondegenerate progression among
the squared differences. For
$\boldsymbol\alpha=(\alpha_1,\ldots,\alpha_L)\in\F_p^L$, define
\[
\begin{aligned}
 \phi_{\boldsymbol\alpha}:\mathbb Z^L&\longrightarrow\F_p,\\
 \mathbf c=(c_1,\ldots,c_L)&\longmapsto
 \sum_{i=1}^L\overline{c_i}\alpha_i,
\end{aligned}
\]
where $\overline{c_i}\in\F_p$ is the residue class of the integer
$c_i$ modulo $p$. Thus, we first reduce the coordinates modulo $p$
and then take their dot product with $\boldsymbol\alpha$ in $\F_p^L$.
We use this map on $\mathcal C$ and on its difference set $\mathcal D$.

It suffices to find a choice of $\boldsymbol\alpha$ for which
we can delete at most
\[
 \frac{|\mathcal C|^2}{p}+\frac{2J(\mathcal C)^3}{p}
\]
points from $\mathcal C$ and satisfy two requirements. First,
$\phi_{\boldsymbol\alpha}$ must be injective on the remaining set.
Second, every three nonzero differences $\delta_1,\delta_2,\delta_3$
of its image must satisfy
\[
 \delta_2^2=\delta_1\delta_3
 \quad\Longrightarrow\quad
 \delta_1^4=\delta_2^4.
\]
By \eqref{eq:GP-difference-criterion}, the the second requirement excludes nondegenerate geometric progressions among the squared differences of that image. Injectivity and the deletion bound give
\eqref{eq:deletion-projection}, so the image will be the required set $A$.

We first identify the quadratic relations that can be left
untreated by the deletion rule.

Write $X=(X_1,\ldots,X_L)$ for the tuple of indeterminates.
For each ordered triple
$\tau=(\mathbf u,\mathbf v,\mathbf w)\in
(\mathcal D\setminus\{\mathbf0\})^3$, define
\[
 P_\tau(X)=(\mathbf v\cdot X)^2
          -(\mathbf u\cdot X)(\mathbf w\cdot X)
 \quad\text{in }\F_p[X_1,\ldots,X_L].
\]
In this polynomial, integer coordinates are reduced modulo $p$
when used as coefficients.
We claim that if $P_\tau\equiv0$, then every projection satisfies
\[
 \phi_{\boldsymbol\alpha}(\mathbf u)^2
 =\phi_{\boldsymbol\alpha}(\mathbf v)^2
 =\phi_{\boldsymbol\alpha}(\mathbf w)^2.
\]
Such an identity cannot yield the nondegenerate relation in
\eqref{eq:GP-difference-criterion}.

Suppose that $P_\tau\equiv0$. Each of the three linear forms
is nonzero: its coefficient vector is a nonzero element of
$\{-3,\ldots,3\}^L$ and $p>18$. The identity
$(\mathbf v\cdot X)^2=(\mathbf u\cdot X)(\mathbf w\cdot X)$
and unique factorization in $\F_p[X_1,\ldots,X_L]$ force
$\mathbf u\cdot X$ and $\mathbf w\cdot X$ to be scalar multiples
of $\mathbf v\cdot X$. We lift this proportionality to $\mathbb Q$
to use the integer digit bounds.
Every $2\times2$ minor formed from their coefficient vectors
vanishes modulo $p$.
For example,
\[
 u_iv_j-u_jv_i\equiv0\pmod p,
 \qquad
 |u_iv_j-u_jv_i|\leq3\cdot3+3\cdot3=18<p.
\]
The same bound holds for the minors formed from any two
of the three vectors. Each minor is divisible by $p$ and has absolute value less than $p$, so it vanishes over $\mathbb Z$.
Thus, the vectors are proportional over $\mathbb Q$, so we may write
\[
 \mathbf u=\lambda_1\mathbf h,\qquad
 \mathbf v=\lambda_2\mathbf h,\qquad
 \mathbf w=\lambda_3\mathbf h,
\]
where $\mathbf h$ is a \emph{primitive} integer vector and
$\lambda_1,\lambda_2,\lambda_3$ are nonzero integers.
Primitivity means $\gcd(h_1,\ldots,h_L)=1$, so every integer
vector on the rational line $\mathbb Q\mathbf h$ is an integer
multiple of $\mathbf h$.
The coordinate bounds give $|\lambda_j|\leq3$, since
$|\lambda_j|\leq|\lambda_jh_i|\leq3$ for any coordinate with $h_i\ne0$.
Since $\mathbf h$ is primitive, it is nonzero modulo $p$,
so $\mathbf h\cdot X$ is a nonzero polynomial. Substitution
into $P_\tau\equiv0$ gives
\[
 (\lambda_2^2-\lambda_1\lambda_3)(\mathbf h\cdot X)^2\equiv0
 \quad\text{in }\F_p[X_1,\ldots,X_L],
\]
and hence
\[
 \lambda_2^2\equiv\lambda_1\lambda_3\pmod p.
\]
Together with
$|\lambda_2^2-\lambda_1\lambda_3|\leq18<p$, this gives the integer
equality $\lambda_2^2=\lambda_1\lambda_3$.
Because $\lambda_2^2\in\{1,4,9\}$ and
$\lambda_1\lambda_3=\lambda_2^2$, the restrictions
$1\leq|\lambda_j|\leq3$ force
$|\lambda_1|=|\lambda_2|=|\lambda_3|$.
For any projection, write
$\delta_1=\phi_{\boldsymbol\alpha}(\mathbf u)$,
$\delta_2=\phi_{\boldsymbol\alpha}(\mathbf v)$, and
$\delta_3=\phi_{\boldsymbol\alpha}(\mathbf w)$.
Linearity gives
\[
 \delta_j=\lambda_j\phi_{\boldsymbol\alpha}(\mathbf h)
 \quad\text{in }\F_p\qquad(j=1,2,3).
\]
Thus, $\delta_1^2=\delta_2^2=\delta_3^2$, proving the claim.

We now choose $\boldsymbol\alpha\in\F_p^L$ uniformly and define
the deletions. For each $\mathbf u\in\mathcal D$, set
\[
 R_{\mathbf u}=\mathcal C\cap(\mathcal C+\mathbf u).
\]
We have $|R_{\mathbf u}|=r_{\mathcal C-\mathcal C}(\mathbf u)$.
Whenever $\mathbf x,\mathbf y\in\mathcal C$ satisfy
$\mathbf x-\mathbf y=\mathbf u$, we have $\mathbf x\in R_{\mathbf u}$, so deleting $R_{\mathbf u}$ leaves no pair with difference
$\mathbf u$.

For each $\tau=(\mathbf u,\mathbf v,\mathbf w)\in
(\mathcal D\setminus\{\mathbf0\})^3$, let
$\mathbf d(\tau)$ be the first entry of \(\tau\) with the smallest difference multiplicity. For a choice of $\boldsymbol\alpha$, define the surviving set $\mathcal C'\subseteq\mathcal C$ by deleting $R_{\mathbf u}$ whenever
$\mathbf u\in\mathcal D\setminus\{\mathbf0\}$ and
$\phi_{\boldsymbol\alpha}(\mathbf u)=0$, and deleting
$R_{\mathbf d(\tau)}$ whenever
$P_\tau\not\equiv0$ and $P_\tau(\boldsymbol\alpha)=0$.

We next estimate the expected number of deleted points.
Every nonzero $\mathbf u\in\mathcal D$ remains nonzero modulo $p$,
since its coordinates have absolute value at most $3$ and $p>18$.
Hence, $\phi_{\boldsymbol\alpha}(\mathbf u)=0$ with probability
$\frac{1}{p}$.
If $P_\tau\not\equiv0$, the degree bound for polynomial zeros
applies with total degree two, giving at most $2p^{L-1}$ zeros
in $\F_p^L$. Thus, $P_\tau(\boldsymbol\alpha)=0$ with probability
at most $2/p$. The number of deleted points is at most the sum of
the cardinalities of the deleted sets, even when these sets overlap.
By linearity of expectation, its expected value is therefore at most
\[
 \frac{1}{p}
 \sum_{\mathbf u\in\mathcal D\setminus\{\mathbf0\}}
r_{\mathcal C-\mathcal C}(\mathbf u)
 +
 \frac{2}{p}
 \sum_{\substack{
       \tau\in(\mathcal D\setminus\{\mathbf0\})^3\\
       P_\tau\not\equiv0}}
r_{\mathcal C-\mathcal C}(\mathbf d(\tau)).
\]
For $\tau=(\mathbf u,\mathbf v,\mathbf w)$, we have
\[
\begin{aligned}
r_{\mathcal C-\mathcal C}(\mathbf d(\tau))
 &=\min\{r_{\mathcal C-\mathcal C}(\mathbf u),
          r_{\mathcal C-\mathcal C}(\mathbf v),
          r_{\mathcal C-\mathcal C}(\mathbf w)\}\\
 &\leq
 \bigl(r_{\mathcal C-\mathcal C}(\mathbf u)
       r_{\mathcal C-\mathcal C}(\mathbf v)
       r_{\mathcal C-\mathcal C}(\mathbf w)\bigr)^{\frac{1}{3}}.
\end{aligned}
\]
The identity
$\sum\limits_{\mathbf u\in\mathcal D\setminus\{\mathbf0\}}
r_{\mathcal C-\mathcal C}(\mathbf u)=|\mathcal C|^2-|\mathcal C|$
counts all ordered pairs of distinct points. Using this identity and enlarging the triple sum to $\mathcal D^3$, we bound the expected number of deletions by
\[
 \frac{|\mathcal C|^2-|\mathcal C|}{p}
 +\frac{2}{p}
  \left(\sum_{\mathbf u\in\mathcal D}
   r_{\mathcal C-\mathcal C}(\mathbf u)^{\frac{1}{3}}
  \right)^3
 =\frac{|\mathcal C|^2-|\mathcal C|+2J(\mathcal C)^3}{p}.
\]
Fix a choice of $\boldsymbol\alpha$ for which the number of deleted
points does not exceed this bound, and set $A=\phi_{\boldsymbol\alpha}(\mathcal C')$.
The restriction of $\phi_{\boldsymbol\alpha}$ to $\mathcal C'$ is
injective: if distinct $\mathbf x,\mathbf y\in\mathcal C'$ had the
same image, then $\mathbf u=\mathbf x-\mathbf y$ would satisfy
$\mathbf u\ne\mathbf0$ and $\phi_{\boldsymbol\alpha}(\mathbf u)=0$, so $\mathbf x\in R_{\mathbf u}$ would have been deleted.
Consequently, $|A|=|\mathcal C'|$. The deletion bound above
is at most $\frac{|\mathcal C|^2+2J(\mathcal C)^3}{p}$, so it proves \eqref{eq:deletion-projection}.

It remains to show that $(A-A)^2$ contains no nondegenerate
geometric progression.
Suppose that nonzero differences
$\delta_1,\delta_2,\delta_3\in A-A$ satisfy
$\delta_2^2=\delta_1\delta_3$.
Choose nonzero vectors
$\mathbf u,\mathbf v,\mathbf w\in\mathcal C'-\mathcal C'$
whose images under $\phi_{\boldsymbol\alpha}$ are
$\delta_1,\delta_2,\delta_3$, respectively.
For $\tau=(\mathbf u,\mathbf v,\mathbf w)$, we then have
$P_\tau(\boldsymbol\alpha)=0$.
If $P_\tau\not\equiv0$, the set $R_{\mathbf d(\tau)}$ would have
been deleted, contradicting
$\mathbf d(\tau)\in\mathcal C'-\mathcal C'$.
Therefore, $P_\tau\equiv0$, and the claim proved above gives
$\delta_1^2=\delta_2^2=\delta_3^2$. In particular,
$\delta_1^4=\delta_2^4$. By \eqref{eq:GP-difference-criterion},
$(A-A)^2$ contains no nondegenerate geometric progression.
This completes the proof.
\end{proof}

The choice of the four digits $0,1,2,3$ ensures that the polynomial
identities left untreated by our deletion rule produce only degenerate progressions. This property fails for five consecutive digits. Indeed, the difference set of $\{0,1,2,3,4\}^L$ contains
$\mathbf h,2\mathbf h,4\mathbf h$, where $\mathbf h$ is a coordinate vector, and
\[
 (2\mathbf h\cdot X)^2-(\mathbf h\cdot X)(4\mathbf h\cdot X)\equiv0.
\]
For any projection with $t=\phi_{\boldsymbol\alpha}(\mathbf h)\ne0$, these differences give the squared values $t^2,4t^2,16t^2$, which form a nondegenerate geometric progression when $p>18$. Since the displayed polynomial is identically zero, it vanishes with probability one, and the $\frac{2}{p}$ bound for nonzero polynomials does not apply. Thus, extending this construction to five consecutive digits would require an additional argument to handle such identities.

\subsection{A family with prescribed digit counts}

To apply Lemma~\ref{lem:digit-projection}, we need a large
family $\mathcal C$ for which $J(\mathcal C)^3$ is small compared
with $p|\mathcal C|$. We prescribe the digit counts so that
the product of the digit weights is constant on the family.
We then bound the weighted pair sums by the corresponding sums
over the full cube, which factor over the coordinates and give
the required estimate for $J(\mathcal C)$.

\begin{proof}[Proof of Proposition~\ref{prop:one-dimensional}\textup{(ii)}]
Let $p$ be a sufficiently large prime. To prove the bound
$|A|\geq p^{\frac{7}{26}}$, we will apply
Lemma~\ref{lem:digit-projection}. For a family
$\mathcal C\subseteq\{0,1,2,3\}^L$, this lemma gives a set
$A\subseteq\F_p$ such that $(A-A)^2$ contains no nondegenerate
geometric progression and
\[
 |A|\geq|\mathcal C|-\frac{|\mathcal C|^2}{p}-\frac{2J(\mathcal C)^3}{p}.
\]

Let $L$ be a positive multiple of $26$, to be chosen below.
For $\mathbf x\in\{0,1,2,3\}^L$, write
$N_j(\mathbf x)=|\{i:x_i=j\}|$ for $j=0,1,2,3$, and set
\begin{equation}\label{eq:fixed-digit-family}
 \mathcal C_L=
 \left\{\mathbf x\in\{0,1,2,3\}^L:
 (N_0(\mathbf x),N_1(\mathbf x),N_2(\mathbf x),N_3(\mathbf x))
 =\frac{L}{26}(5,8,8,5)\right\}.
\end{equation}
We will apply the lemma with $\mathcal C=\mathcal C_L$ and
choose $L$ so that
\begin{equation}\label{eq:fixed-digit-targets}
 |A|\geq\frac{1}{2}|\mathcal C_L|,
 \qquad
 \frac{1}{2}|\mathcal C_L|\geq p^{\frac{7}{26}}.
\end{equation}
For the first inequality, it suffices to prove
\[
 \frac{|\mathcal C_L|^2}{p}+\frac{2J(\mathcal C_L)^3}{p}
 \leq\frac{1}{2}|\mathcal C_L|.
\]
We therefore begin with estimates for $|\mathcal C_L|$ and
$J(\mathcal C_L)$.

The prescribed digit counts give
\[
 |\mathcal C_L|
 =\frac{L!}
 {\left(\left(\frac{5L}{26}\right)!\right)^2
  \left(\left(\frac{8L}{26}\right)!\right)^2}.
\]
Applying Stirling's formula in the form
\[
 \log(n!)=n\log n-n+\frac{1}{2}\log n+O(1),
\]
we obtain
\[
 \log|\mathcal C_L|
 =h_0L-\frac{3}{2}\log L+O(1),
 \qquad
 h_0=\log26-\frac{5}{13}\log5-\frac{8}{13}\log8.
\]
In particular,
\begin{equation}\label{eq:fixed-digit-size}
 |\mathcal C_L|\asymp L^{-\frac{3}{2}}e^{h_0L},
\end{equation}
with absolute implied constants.

We next estimate $J(\mathcal C_L)$. Set
$(w_0,w_1,w_2,w_3)=(1,3,3,1)$ and define
\[
 \omega(\mathbf x)=\prod_{i=1}^L w_{x_i},
 \qquad \mathbf x\in\{0,1,2,3\}^L.
\]
The prescribed counts make this weight constant on $\mathcal C_L$:
\[
 \omega(\mathbf x)
 =3^{N_1(\mathbf x)+N_2(\mathbf x)}
 =3^{\frac{8L}{13}}
 \qquad(\mathbf x\in\mathcal C_L).
\]
For $j\in\{-3,\ldots,3\}$, let
\[
 \rho(j)=\sum_{\substack{a,b\in\{0,1,2,3\}\\a-b=j}}w_aw_b.
\]
These are the weighted pair counts for one coordinate. Direct
summation gives
\begin{equation}\label{eq:fixed-digit-weight-table}
 \rho(0)=20,\qquad
 \rho(\pm1)=15,\qquad
 \rho(\pm2)=6,\qquad
 \rho(\pm3)=1.
\end{equation}

For each $\mathbf u\in\mathcal C_L-\mathcal C_L$, the difference
multiplicity defined in Lemma~\ref{lem:digit-projection} satisfies
\begin{equation}\label{eq:fixed-digit-weight-bound}
\begin{aligned}
 3^{\frac{16L}{13}}
 r_{\mathcal C_L-\mathcal C_L}(\mathbf u)
 &=\sum_{\substack{\mathbf x,\mathbf y\in\mathcal C_L\\
                   \mathbf x-\mathbf y=\mathbf u}}
       \omega(\mathbf x)\omega(\mathbf y)\\
 &\leq
 \sum_{\substack{\mathbf x,\mathbf y\in\{0,1,2,3\}^L\\
                 \mathbf x-\mathbf y=\mathbf u}}
       \omega(\mathbf x)\omega(\mathbf y)\\
 &=\prod_{i=1}^L\rho(u_i).
\end{aligned}
\end{equation}
The inequality enlarges the sum to the full cube, where the
coordinates can be summed independently.

Taking cube roots and enlarging the sum over differences to
$\{-3,\ldots,3\}^L$, we obtain
\[
\begin{aligned}
 J(\mathcal C_L)
 &=\sum_{\mathbf u\in\mathcal C_L-\mathcal C_L}
       r_{\mathcal C_L-\mathcal C_L}(\mathbf u)^{\frac{1}{3}}\\
 &\leq 3^{-\frac{16L}{39}}
       \sum_{\mathbf u\in\{-3,\ldots,3\}^L}
       \prod_{i=1}^L\rho(u_i)^{\frac{1}{3}}\\
 &=3^{-\frac{16L}{39}}
   \left(\sum_{j=-3}^3\rho(j)^{\frac{1}{3}}\right)^L
 =3^{-\frac{16L}{39}}s_0^L,
\end{aligned}
\]
where
\[
 s_0:=\sum_{j=-3}^3\rho(j)^{\frac{1}{3}}
 =\sqrt[3]{20}+2\sqrt[3]{15}+2\sqrt[3]{6}+2.
\]
Cubing the preceding estimate gives
\begin{equation}\label{eq:fixed-digit-moment}
 J(\mathcal C_L)^3
 \leq\left(3^{-\frac{16}{13}}s_0^3\right)^L
 =e^{b_0L},
 \qquad
 b_0:=3\log s_0-\frac{16}{13}\log3.
\end{equation}
Since $s_0>12$ and $0<h_0\leq\log4$, we have
\[
 b_0>3\log12-\frac{16}{13}\log3
 >2\log4\geq2h_0.
\]

We now choose $L$ to prove the first inequality in
\eqref{eq:fixed-digit-targets}. Dividing the two error terms by
$|\mathcal C_L|$ and using \eqref{eq:fixed-digit-size} and
\eqref{eq:fixed-digit-moment}, we obtain
\[
\begin{aligned}
 \frac{|\mathcal C_L|}{p}
 +\frac{2J(\mathcal C_L)^3}{p|\mathcal C_L|}
 &\ll
 \frac{e^{h_0L}}{pL^{\frac{3}{2}}}
 +\frac{L^{\frac{3}{2}}e^{(b_0-h_0)L}}{p}\\
 &\ll
 \frac{L^{\frac{3}{2}}e^{(b_0-h_0)L}}{p},
\end{aligned}
\]
where the last step use $b_0>2h_0$ and $L\geq1$.

For $L\leq\frac{\log p}{b_0-h_0}$, we have
$L^{\frac{3}{2}}\ll(\log p)^{\frac{3}{2}}$.
It therefore suffices to require
\[
 e^{(b_0-h_0)L}
 \leq e^{-K}p(\log p)^{-\frac{3}{2}},
\]
where $K>0$ will be fixed sufficiently large.
Taking logarithms, we choose the largest multiple of $26$
satisfying this inequality, namely
\begin{equation}\label{eq:fixed-digit-length}
 L=26\left\lfloor
 \frac{\log p-\frac{3}{2}\log\log p-K}
      {26(b_0-h_0)}
 \right\rfloor.
\end{equation}
For $p$ sufficiently large in terms of $K$, this choice is
positive, satisfies $L\leq\frac{\log p}{b_0-h_0}$, and has
$L\asymp\log p$. Consequently,
\[
 \frac{|\mathcal C_L|}{p}
 +\frac{2J(\mathcal C_L)^3}{p|\mathcal C_L|}
 \ll
 \frac{L^{\frac{3}{2}}e^{(b_0-h_0)L}}{p}
 \leq
 \frac{e^{-K}}{(b_0-h_0)^{\frac{3}{2}}},
\]
with an implied constant independent of $K$.
We first fix $K$ sufficiently large and then take $p$
sufficiently large so that the relative error is at most
$\frac{1}{2}$. Lemma~\ref{lem:digit-projection} then gives
\[
 |A|\geq\frac{1}{2}|\mathcal C_L|.
\]
This proves the first inequality in
\eqref{eq:fixed-digit-targets}.

To prove the second inequality, put
\[
 \beta_0:=\frac{h_0}{b_0-h_0} \in (0,1).
\]
With $K$ now fixed, rounding in \eqref{eq:fixed-digit-length}
changes $L$ by less than $26$, and
$\log L=\log\log p+O(1)$.
Substituting this choice of $L$ into the logarithmic size
estimate gives
\[
 \log|\mathcal C_L|
 =\beta_0\log p
  -\frac{3}{2}(1+\beta_0)\log\log p+O(1).
\]
Hence,
\begin{equation}\label{eq:fixed-digit-lower-bound}
 \frac{1}{2}|\mathcal C_L|
 \gg
 \frac{p^{\beta_0}}
      {(\log p)^{\frac{3}{2}(1+\beta_0)}}.
\end{equation}

It remains to verify that $\beta_0>\frac{7}{26}$.
Since $b_0-h_0>0$, this is equivalent to $33h_0>7b_0$.
Substituting the definitions yields
\begin{equation}\label{eq:fixed-digit-exponent-check}
 13(33h_0-7b_0)
 =\log\left(
 \frac{26^{429}\cdot3^{112}}
      {5^{165}\cdot8^{264}\cdot s_0^{273}}
 \right).
\end{equation}
We verify this inequality using rational upper bounds for the
cube roots and integer arithmetic. The inequalities
\[
 2715^3>20\cdot1000^3,\qquad
 2467^3>15\cdot1000^3,\qquad
 1818^3>6\cdot1000^3
\]
imply $s_0<\frac{13285}{1000}$, while direct integer
arithmetic gives
\[
 26^{429}\cdot3^{112}\cdot1000^{273}
 >
 5^{165}\cdot8^{264}\cdot13285^{273}.
\]
Thus the ratio in \eqref{eq:fixed-digit-exponent-check}
is greater than $1$, proving $\beta_0>\frac{7}{26}$.

Since
\[
 \frac{p^{\beta_0-\frac{7}{26}}}
      {(\log p)^{\frac{3}{2}(1+\beta_0)}}
 \longrightarrow\infty
 \qquad(p\to\infty),
\]
the estimate \eqref{eq:fixed-digit-lower-bound} implies
$\frac{1}{2}|\mathcal C_L|\geq p^{\frac{7}{26}}$
for all sufficiently large $p$.
This proves the second inequality in
\eqref{eq:fixed-digit-targets}, and the first now gives
$|A|\geq p^{\frac{7}{26}}$.

The same set $A$ works for every $\kappa\in\F_p^\times$:
if $\kappa(A-A)^2$ contained a nondegenerate geometric
progression, multiplication by $\kappa^{-1}$ would give one
in $(A-A)^2$. This completes the proof.
\end{proof}

Numerically, $\beta_0\approx0.2693287509$, so $\frac{7}{26}$
is a a convenient rational exponent below $\beta_0$.
We do not claim that these choices or the resulting exponent
are optimal.

\section{Further remarks}\label{sec:remarks}

The preceding sections reduce the odd-dimensional extremal
problem to dimension one. We now explain what is known in even
dimensions, where the split and nonsplit forms behave differently.
We then show why the condition $|D|\geq p^{1-\varepsilon}$, for fixed $\varepsilon>0$, does not force a progression in an arbitrary set $D\subset\F_p$, and state the remaining questions.

\subsection{Even dimensions}\label{sec:even-remarks}

The distance estimate of Pham, Shen, Tran, and Xue
\cite{PhamShenTranXue2026}, combined with a quantitative form of
Roth's theorem, gives a uniform sufficient threshold in even
dimensions. This threshold is sharp for split forms. For the
standard nonsplit form, the remaining extremal problem reduces
to the plane.
We first bound the cardinality of an arbitrary set containing
no nondegenerate geometric progression. This will be applied
to distance sets in both cases.

\Needspace{14\baselineskip}
\begin{lemma}\label{lem:positive-density-gp}
There is an absolute constant $a>0$ such that, for every
sufficiently large prime $p$, every set $D\subset\F_p$ containing
no nondegenerate geometric progression satisfies
\begin{equation}\label{eq:gp-free-cardinality}
 |D|\leq p \, e^{-a\left(\frac{\log p}{\log\log p}\right)^{\frac{1}{6}}}.
\end{equation}
Consequently, for every $c,K>0$, there exists $p_0=p_0(c,K)$
such that
\[
 |D|\geq\frac{cp}{(\log p)^K}
 \quad\Longrightarrow\quad
 D\text{ contains a nondegenerate geometric progression}
\]
whenever $p\geq p_0$ is prime and $D\subset\F_p$.
\end{lemma}
\begin{proof}
Let $r_3(N)$ be the largest cardinality of a subset of
$\{0,\ldots,N-1\}$ containing no three distinct terms in arithmetic
progression. By Raghavan~\cite[Theorem~1.4]{Raghavan2026},
there is an absolute constant $b>0$ such that
\[
 r_3(N)\leq N \, e^{-b\left(\frac{\log N}{\log\log N}\right)^{\frac{1}{6}}}
\]
for all sufficiently large $N$.

Suppose that $D$ contains no nondegenerate geometric progression.
Choose a primitive root $\gamma\in\F_p^\times$ and let
\[
 S=\{j\in\{0,\ldots,p-2\}:\gamma^j\in D\}.
\]
If $S$ contained integers $u,u+v,u+2v$ with $v>0$, then
$\gamma^u,\gamma^{u+v},\gamma^{u+2v}$ would be three distinct nonzero elements
of $D$ in geometric progression. Their distinctness follows
because their exponents belong to $\{0,\ldots,p-2\}$.
Thus, $|S|\leq r_3(p-1)$ and
\[
 |D|\leq1+r_3(p-1)
 \leq p \, e^{-a\left(\frac{\log p}{\log\log p}\right)^{\frac{1}{6}}}
\]
for all sufficiently large $p$, with a suitable absolute
constant $a>0$.
The final assertion follows since
$e^{-a\left(\frac{\log p}{\log\log p}\right)^{\frac{1}{6}}}=o((\log p)^{-K})$ for every fixed
$K>0$. This completes the proof.
\end{proof}

\paragraph{A uniform sufficient threshold.}
We apply the preceding lemma to the recent distance estimate.
For a quadratic form $Q$ on $\F_p^{2m}$, let
\[
 \Delta_Q(E)=\{Q(x-y):x,y\in E\}.
\]
A nondegenerate form on $\F_p^{2m}$ is called \emph{split} if
it has an $m$-dimensional totally isotropic subspace, and
\emph{nonsplit} otherwise. In the nonsplit case, the largest
dimension of a totally isotropic subspace is $m-1$.
Pham, Shen, Tran, and Xue~\cite[Theorem~1.1 and Remark~1.2]{PhamShenTranXue2026}
proved that, for every nondegenerate $Q$ and every
$E\subset\F_p^{2m}$ with $|E|>p^m$,
\begin{equation}\label{eq:even-recent-distance-bound}
 |\Delta_Q(E)|\geq\frac{c_0p}{\log^2(2p)},
\end{equation}
where $c_0>0$ is absolute. If $Q$ is nonsplit, the same conclusion
holds when $|E|=p^m$; in fact, their estimate gives
$|\Delta_Q(E)|\geq \frac{c_0 p}{\log(2p)}$ in this case, after decreasing
$c_0$ if necessary. The constants are independent of $m$.

\begin{corollary}\label{cor:even-dimensional-gp}
There is an absolute constant $p_0$ such that the following holds
for every odd prime $p\geq p_0$ and every integer $m\geq1$.
Let $Q$ be a nondegenerate quadratic form on $\F_p^{2m}$ and let
$E\subset\F_p^{2m}$. If $|E|>p^m$, then $\Delta_Q(E)$ contains a
nondegenerate geometric progression. The same conclusion holds
if $|E|=p^m$ and $Q$ is nonsplit.
\end{corollary}
\begin{proof}
For $p\geq2$, we have $\log(2p)\leq2\log p$. Thus,
\eqref{eq:even-recent-distance-bound} gives
\[
 |\Delta_Q(E)|\geq\frac{c_0p}{4(\log p)^2}.
\]
Lemma~\ref{lem:positive-density-gp}, with $c= \frac{1}{4}c_0$ and
$K=2$, proves the assertion when $|E|>p^m$. The quoted estimate
for nonsplit forms gives the conclusion at equality as well.
Since $c_0$ is absolute, the lower bound on $p$ is independent of
$m$ and $Q$. This completes the proof.
\end{proof}

\paragraph{Split forms.}
The strict inequality $|E|>p^m$ is sharp when $Q$ is split.
Indeed, a split form has a totally isotropic subspace $W$ of
dimension $m$. Taking $E=W$ gives $|E|=p^m$ and
$\Delta_Q(E)=\{0\}$. For the standard form $Q_{2m}$, the split
condition is that $(-1)^m$ be a square in $\F_p$; see
\cite[Lemma~2.7]{PhamShenXue2026}. Consequently,
\begin{equation}\label{eq:even-split-exact}
 M_{\GP,2m}(p)=p^m
\end{equation}
for all sufficiently large primes $p$ whenever $m$ is even or
$p\equiv1\pmod4$.

\paragraph{Nonsplit forms.}
For the standard form, the remaining case is $p\equiv3\pmod4$
and $m$ odd. A totally isotropic subspace then has dimension at
most $m-1$, so the preceding example does not establish
sharpness of the exponent $m$.
A cylinder over a planar set gives a lower comparison with
$M_{\GP,2}(p)$, and the extraction theorem gives the upper
comparison. We record both inequalities together with the
available planar bounds.

\begin{proposition}\label{prop:even-nonsplit-comparison}
There are absolute constants $c>0$, $C\geq1$, and $p_0$ such
that, for every prime $p\geq p_0$ with $p\equiv3\pmod4$ and
every odd integer $m\geq1$,
\begin{equation}\label{eq:even-nonsplit-comparison}
 p^{m-1}M_{\GP,2}(p)
 \leq M_{\GP,2m}(p)
 \leq Cp^{m-1}M_{\GP,2}(p).
\end{equation}
Moreover,
\begin{equation}\label{eq:planar-gp-bounds}
 p^{\frac{7}{26}}\leq M_{\GP,2}(p)
 \leq Cp \, e^{-c\left(\frac{\log p}{\log\log p}\right)^{\frac{1}{6}}}.
\end{equation}
\end{proposition}
\begin{proof}
We first prove \eqref{eq:even-nonsplit-comparison} and then
establish the planar bounds in \eqref{eq:planar-gp-bounds}.

By the classification of quadratic forms over $\F_p$, a linear
change of coordinates identifies $Q_{2m}$ with
\[
 u\cdot v+Q_2(z),
 \qquad u,v\in\F_p^{m-1},\quad z\in\F_p^2;
\]
see \cite[Lemma~2.7]{PhamShenXue2026}. If $P\subset\F_p^2$ has
no nondegenerate geometric progression in its distance set,
then, in these coordinates, the set
\[
 \F_p^{m-1}\times\{0\}^{m-1}\times P
\]
has $p^{m-1}|P|$ elements and distance set $\Delta(P)$.
This proves the lower bound in
\eqref{eq:even-nonsplit-comparison}.

For the upper bound, let $E\subset\F_p^{2m}$ have no
nondegenerate geometric progression in $\Delta(E)$.
Corollary~\ref{cor:even-dimensional-gp} gives $|E|<p^m$.
The extraction theorem of Pham, Shen, and Xue
\cite[Theorem~1.1]{PhamShenXue2026}, with the residual binary
form identified with $Q_2$, gives a set $P\subset\F_p^2$ such
that $\Delta(P)\subseteq\Delta(E)$ and
\[
 |P|\geq
 \frac{1}{C_0}\frac{|E|/p^{m-1}}{1+|E|/p^{m+1}}
 \geq\frac{|E|}{2C_0p^{m-1}},
\]
where $C_0$ is absolute. Since $|P|\leq M_{\GP,2}(p)$,
the upper bound follows.

To prove \eqref{eq:planar-gp-bounds}, let $A\subset\F_p$ be
the set from Proposition~\ref{prop:one-dimensional}\textup{(ii)}. Then
$\Delta(A\times\{0\})=(A-A)^2$, giving
$M_{\GP,2}(p)\geq p^{\frac{7}{26}}$.
Next, let $P\subset\F_p^2$ be nonempty and have no nondegenerate geometric
progression in $\Delta(P)$. Again,
Corollary~\ref{cor:even-dimensional-gp} gives $|P|<p$.
Since $Q_2$ is anisotropic, the planar estimate of
Pham, Shen, Tran, and Xue
\cite[Theorem~1.3]{PhamShenTranXue2026} yields
\[
 \frac{|P|}{\log(2|P|)}\ll|\Delta(P)|
 \leq p \, e^{-a\left(\frac{\log p}{\log\log p}\right)^{\frac{1}{6}}},
\]
where the last inequality is
Lemma~\ref{lem:positive-density-gp}.
Since $|P|<p$, we obtain
\[
\begin{aligned}
 |P|&\ll p\log(2p) \, e^{-a\left(\frac{\log p}{\log\log p}\right)^{\frac{1}{6}}}
 \ll p \, e^{-c\left(\frac{\log p}{\log\log p}\right)^{\frac{1}{6}}}
\end{aligned}
\]
for a suitable absolute constant $c>0$.
This completes the proof.
\end{proof}

In particular, for $p\equiv3\pmod4$ and odd $m$,
\begin{equation}\label{eq:even-nonsplit-bounds}
 p^{m-1+\frac{7}{26}}\leq M_{\GP,2m}(p)
 \ll p^m \, e^{-c\left(\frac{\log p}{\log\log p}\right)^{\frac{1}{6}} }.
\end{equation}
Thus, $M_{\GP,2m}(p)=o(p^m)$ as $p$ tends to infinity through
these primes. The upper bound does not give a fixed power
improvement below $p^m$, and the lower bound does not establish
sharpness of the exponent $m$. Determining the optimal exponent
in this nonsplit case reduces, by
\eqref{eq:even-nonsplit-comparison}, to determining it for
$M_{\GP,2}(p)$ when $p\equiv3\pmod4$.

\subsection{Cardinality alone is not enough}\label{sec:cardinality-alone}

Lemma~\ref{lem:positive-density-gp} forces a progression at
every density $c/(\log p)^K$. A fixed power bound is different:
for every $\varepsilon>0$, the size condition
$|D|\geq p^{1-\varepsilon}$ does not force a nondegenerate
geometric progression in an arbitrary set $D\subset\F_p$.
Indeed, let $\ell_0=\left\lfloor\frac{p-1}{3}\right\rfloor$.
Behrend's construction \cite{Behrend1946} gives a set
$B\subset\{0,\ldots,\ell_0-1\}$ with
$|B|\geq \ell_0e^{-C\sqrt{\log \ell_0}}$ containing no arithmetic progression of three distinct terms.
Choose a generator $\gamma$ of $\F_p^\times$ and set
$D=\{\gamma^b:b\in B\}$.
A nondegenerate geometric progression in $D$ would give distinct
$b_1,b_2,b_3\in B$ satisfying $b_1+b_3\equiv2b_2\pmod{p-1}$.
Since $|b_1+b_3-2b_2|\leq2(\ell_0-1)<p-1$, this congruence would imply $b_1+b_3=2b_2$ over the integers, contrary to the choice of $B$.

After increasing the absolute constant $C$ if necessary, we obtain
\[
 |D|=|B|\geq p e^{-C\sqrt{\log p}}
\]
for all sufficiently large $p$.
For every fixed $\varepsilon>0$, this lower bound exceeds
$p^{1-\varepsilon}$ once $p$ is sufficiently large, showing that
such a size condition alone does not force a nondegenerate
geometric progression.
The argument in Section~\ref{sec:one-dimensional-upper}
uses the additional structure of squared difference sets.

\subsection{Remarks over arbitrary finite fields}
Let $q$ be a power of a prime $p>18$. Define $\Delta(E)$ and
$M_{\GP,d}(q)$ over $\F_q$ by the same formulas as over $\F_p$.
The arguments in Sections~\ref{sec:reduction}--\ref{sec:odd-examples}
extend to this setting. More precisely, the quadratic-form counts and the upper-bound argument use $q$ in place
of $p$, with Fourier analysis defined using a nontrivial additive character
of $\F_q$. In Lemma~\ref{lem:digit-projection}, choosing
$\boldsymbol\alpha$ uniformly from $\F_q^L$ gives the probabilities
$1/q$ for a collision and at most $2/q$ for a nonzero quadratic relation.
The two steps that lift congruences to integer equalities still use
the characteristic bound $p>18$.
The digit-count estimates are unchanged, with $L$ chosen in terms of $\log q$.
For every integer $m\geq0$,
\[
 q^mM_{\GP,1}(q)\leq M_{\GP,2m+1}(q)
 \leq4q^mM_{\GP,1}(q).
\]
Moreover, there are absolute constants $C\geq1$ and $q_0$ such that
\[
 q^{m+\frac7{26}}\leq M_{\GP,2m+1}(q)
 \leq Cq^{m+\frac{16}{25}}
\]
whenever $q\geq q_0$. These constants are independent of $m$,
$p$, and the degree of the field extension.

Questions of sharpness and the construction of extremal examples are
more involved over general finite fields. To keep the presentation
focused and consistent with recent progress on the distance problem
over prime fields, we state and prove the main results only over prime
fields.

\subsection{Open questions}\label{sec:remaining-questions}

\paragraph{The one-dimensional extremal problem.}
Proposition~\ref{prop:one-dimensional} gives, for all sufficiently
large primes $p$,
\[
 p^{\frac{7}{26}}\leq M_{\GP,1}(p)\ll p^{\frac{16}{25}}.
\]
The gap between these bounds leads to the following question.

\begin{question}\label{question:one-dimensional-extremal}
What is the order of magnitude of $M_{\GP,1}(p)$ as $p$ tends
to infinity through the primes?
\end{question}

By Proposition~\ref{prop:odd-dimension-comparison},
\[
 M_{\GP,2m+1}(p)\asymp p^mM_{\GP,1}(p),
\]
with absolute implied constants, so an answer to
Question~\ref{question:one-dimensional-extremal} would also
determine the extremal order in every odd dimension.

\paragraph{Other arithmetic patterns.}
A natural next question is how large $E$ must be for $\Delta(E)$ to contain other patterns, such as a nontrivial three-term arithmetic progression or a geometric progression of length at least four.

\paragraph{On $(u,s)$--Salem sets.}
Fraser introduced $(u,s)$--Salem sets and studied their distance sets in \cite{fraser2}; see also \cite{CGKPTZ25} for subsequent results. It would be interesting to determine whether suitable $(u,s)$--Salem hypotheses lower the size thresholds for geometric progressions proved in this paper.

{\bf Acknowledgement.} 
D. T. Tran would like to thank the Vietnam Institute for Advanced Study in Mathematics (VIASM) for its warm hospitality and excellent research environment.

\end{document}